\documentclass[final,fleqn,3p,12pt,times]{elsarticle}
\usepackage[T1]{fontenc}
\usepackage{amsmath,mathdots}

\biboptions{sort&compress}

\usepackage{amsmath}
\makeatletter
\let\oldtheequation\theequation
\renewcommand\tagform@[1]{\maketag@@@{\ignorespaces#1\unskip\@@italiccorr}}
\renewcommand\theequation{(\oldtheequation)}
\makeatother
\usepackage{amsthm}
\usepackage[below]{placeins}
\usepackage{paralist}

\usepackage{tabularx}
\usepackage{multirow}

\usepackage[defaultcolor=red]{changes}

\usepackage{aliascnt}

\newproof{pf}{Proof}

\newaliascnt{lemma}{theorem}

\aliascntresetthe{lemma}

\newaliascnt{corollary}{theorem}

\aliascntresetthe{corollary}

\DeclareMathOperator{\diag}{diag}

\newcommand{\coltext}{(For interpretation of the references to color in this figure legend, the reader is referred to the web version of this article.)}

\usepackage{hyperref}
\hypersetup{
        pdfinfo = {
        Title =    {Integrating adaptive optimization into least squares progressive iterative approximation},
        Authors =  {},
        Subject =  {},
        Keywords = {},
        Creator =  {},
        Producer = {}
    },
    pdfstartview =  FitH,
    bookmarksopen = True,
    colorlinks =    True
}

\makeatletter
\long\def\@makecaption#1#2{%
   \vskip\abovecaptionskip
   \sbox\@tempboxa{\footnotesize\textbf{#1.} #2}%
   \ifdim \wd\@tempboxa >\hsize
     \footnotesize{\textbf{#1.} #2}\par
   \else
     \global \@minipagefalse
     \hb@xt@\hsize{\hfil\box\@tempboxa\hfil}%
   \fi
   \vskip\belowcaptionskip}
\makeatother

\begin{document}

\begin{frontmatter}

\title{Integrating adaptive optimization into least squares progressive iterative approximation}

\begin{abstract}
    This paper introduces the Adaptive Gradient Least Squares Progressive iterative Approximation (AdagradLSPIA), an accelerated version of the Least Squares Progressive Iterative Approximation (LSPIA) method, enhanced with adaptive optimization techniques inspired by the adaptive gradient (Adagrad) algorithm. By using historical (accumulated) gradient information to dynamically adjust weights, AdagradLSPIA achieves faster convergence compared to the standard LSPIA method. The effectiveness of AdagradLSPIA is demonstrated through its application to tensor product B-spline surface fitting, where this method consistently outperforms LSPIA in terms of accuracy, computational efficiency, and robustness to variations in global weight selection.
\end{abstract}


\author{Svaj\={u}nas Sajavi\v{c}ius}
\address{
  Department of Software Engineering, Faculty of Informatics,
  Kaunas University of Technology, Lithuania
}
\ead{svajunas.sajavicius@ktu.lt}
\ead[url]{https://sajavicius.github.io/}

\end{frontmatter}

\section{Introduction}
\label{sec:1}

In geometric design, data fitting, and reverse engineering, creating accurate and efficient representations of complex shapes and surfaces from given data points is a fundamental task. This is particularly important in applications such as computer graphics, CAD/CAM systems, and digital manufacturing, where precise surface fitting plays a key role in ensuring high-quality modeling, simulation, and fabrication. Among the various methods developed for these purposes, Least Squares Progressive Iterative Approximation (LSPIA) has attracted significant attention due to its simplicity and effectiveness in handling large data sets.

LSPIA, initially introduced by~\citet{Deng2014}, is a geometric iterative method (GIM) that extends the standard Progressive Iterative Approximation (PIA)~\cite{Lin2004} by minimizing the least squares error between the data points and the fitting patch. The method constructs a sequence of fitting patches by iteratively updating the control points, thereby reducing the fitting error at each step. This approach has been widely adopted for various geometric modeling tasks, including mesh generation, solid modeling, and reverse engineering~\cite{Lin2018}. Despite its effectiveness, the standard LSPIA method relies on a constant global weight to scale adjusting vectors, which can limit the convergence speed and accuracy, especially when dealing with data sets that exhibit significant variability.

Recent papers have introduced several enhancements to the LSPIA method, focusing on improvements in convergence speed and computational efficiency. For example, memory-augmented methods such as MLSPIA~\cite{Huang2020} have been proposed, demonstrating superior performance through advanced weighting mechanisms and relaxation strategies. The asynchronous LSPIA (ALSPIA)~\cite{Wu2024} method significantly reduces iteration counts and computational time compared to the standard LSPIA, particularly in large-scale problems, by leveraging asynchronous updates to the adjusting vectors. Another notable improvement is the introduction of momentum techniques in LSPIA for curve fitting problem, as seen in the PmLSPIA and NmLSPIA methods~\cite{Liu2024}, which utilize Polyak's and Nesterov's momentum to achieve faster convergence by better determining the search direction. Such techniques as parameter correction~\cite{Song2024}, high-order matrix inverse approximation~\cite{Ebrahimi2019,Sajavicius2023,Yao2025}, matrix distributing~\cite{Yao2024}, or NURBS weights and knots optimization~\cite{Lan2024} have also been employed to further optimize the convergence rates and computational efficiency of LSPIA. These advancements highlight the continuous evolution of LSPIA and its adaptability to address complex geometric modeling tasks. 

One way to improve the performance of GIMs like LSPIA is by adopting adaptive optimization techniques. Adaptive optimization methods, such as the Adaptive Gradient (Adagrad) algorithm~\cite{Duchi2011}, have become popular in machine learning due to their ability to dynamically adjust learning rates based on historical (accumulated) gradient information (see also~\cite{Bottou2018,Kochenderfer2019}). This adaptability allows for faster and more stable convergence, particularly when dealing with data where gradient magnitudes vary significantly across parameters. Motivated by this approach, an accelerated version of the standard LSPIA, referred to as Adaptive Gradient Least Squares Progressive Iterative Approximation (AdagradLSPIA), is proposed in this paper.

The AdagradLSPIA method integrates the adaptive weighting mechanism of Adagrad with the geometric iterative framework of LSPIA. Unlike the standard LSPIA approach, where a fixed global weight is used, AdagradLSPIA assigns separate adaptive weights to each adjusting vector based on accumulated gradient information from previous iterations. This not only improves the convergence rate by allowing larger weights in smoother regions of the surface, but also provides more stability when handling complex geometric features and noisy data.

The remainder of this paper is organized as follows. \autoref{sec:2}~presents the formulation of the AdagradLSPIA method, including a detailed discussion of the underlying geometric iterative procedure and the adaptive weighting mechanism, as well as remarks on convergence, hyperparameter tuning, and computational complexity. This is followed by an extensive experimental study in~\autoref{sec:3} that demonstrates the effectiveness of the proposed method and its
suitability for surface fitting tasks in geometric modeling and reverse engineering. Finally, \autoref{sec:4}~concludes the paper with a summary of the findings and suggestions for future research.

\section{Adaptive gradient least squares progressive iterative approximation}
\label{sec:2}

In this section, the Adaptive Gradient Least Squares Progressive Iterative Approximation (AdagradLSPIA) method is introduced, which integrates adaptive optimization techniques into the traditional LSPIA framework. The data fitting problem is first formulated, and the geometric iterative method for its solution is outlined. A detailed explanation of the adaptive weighting mechanism, which drives the convergence of the method, is then provided. Key aspects of convergence, hyperparameter tuning, and computational complexity are also discussed.

\subsection{Least squares fitting problem}
The goal of the AdagradLSPIA method is to fit a smooth curve or surface to a given set of data points by minimizing the least squares error.

Let $\mathcal{I}$ and $\mathcal{J}$ be appropriately chosen (multi-)index sets such that $\lvert \mathcal{I} \rvert \leq \lvert \mathcal{J} \rvert$. Consider a set of data points $\{\mathbf{q}_j\}_{j \in \mathcal{J}} \subset \mathbb{R}^d$ (where $d = 2, 3$) with corresponding parameter values $\{\mathbf{t}_j\}_{j \in \mathcal{J}}$ in a suitable parameter space $\Omega \subset \mathbb{R}$ (or $\mathbb{R}^2$ for surfaces). The \emph{fitting patch}, i.e.~a curve or surface that approximates the given data points, is represented as
\begin{equation*}
    \mathbf{c}(\mathbf{t}) = \sum_{i \in \mathcal{I}}{\mathbf{p}_i B_i(\mathbf{t})}, \quad \mathbf{t} \in \Omega,
\end{equation*}
where $\mathbf{p}_i \in \mathbb{R}^d$ are the \emph{control points} to be determined, and $B_i(\mathbf{t})$ are \emph{blending functions} (e.g., B-spline or NURBS basis functions).

The considered fitting problem can be formulated as the least squares problem
\begin{equation}
\label{eq:obj}
    E(\mathbf{p}) = \frac{1}{2}\sum_{j \in \mathcal{J}}{\lVert \mathbf{q}_j - \mathbf{c}(\mathbf{t}_j) \rVert_2^2} \rightarrow \min
\end{equation}
where $E(\mathbf{p})$ is the \emph{objective function} (least squares error), and $\mathbf{p} = [\mathbf{p}_i]_{i \in \mathcal{I}}^T$ is the vector containing the control points of the fitting patch $\mathbf{c}(\mathbf{t})$.

The matrix form of the objective function is
\begin{equation*}
    E(\mathbf{p}) = \frac{1}{2} \lVert \mathbf{A} \mathbf{p} - \mathbf{q} \rVert_2^2
\end{equation*}
where $\mathbf{A} = (\mathbf{A}_{ij})$ is the collocation matrix containing basis function values, $\mathbf{A}_{ij} = B_j(\mathbf{t}_i)$, and $\mathbf{q} = [\mathbf{q}_j]_{j \in \mathcal{J}}^T$ is the vector of given data points. The gradient of $ E(\mathbf{p})$ is given by
\begin{equation*}
    \nabla E(\mathbf{p}) = \mathbf{A}^T(\mathbf{A} \mathbf{p} - \mathbf{q}).
\end{equation*}

Since $E(\mathbf{p})$ is a quadratic function of $\mathbf{p}$, it is always convex. For $E(\mathbf{p})$ to be strictly convex, the Hessian
\begin{equation*}
    \nabla^2 E(\mathbf{p}) = \mathbf{A}^T \mathbf{A}
\end{equation*}
must be positive definite, which happens if and only if  the Gram matrix $\mathbf{A}^T \mathbf{A}$ has full column rank (i.e., the blending functions used in the collocation are linearly independent).

B-spline and NURBS basis functions are typically locally supported, which means that each column of $\mathbf{A}$ represents a weighted contribution of basis functions at different collocation points. If the control points are sufficient and properly distributed, and the basis functions are linearly independent, then $\mathbf{A}$ has the full column rank, which implies that $\mathbf{A}^T \mathbf{A}$ is positive definite and $E(\mathbf{p})$ is strictly convex. In this case, there exists a unique minimizer for $\mathbf{p}$.

If the degree of the basis functions is too low or the knots are poorly chosen, some basis functions may be redundant or nearly dependent, leading to a rank-deficient $\mathbf{A}$. If $\mathbf{A}$ does not have full column rank, then $\mathbf{A}^T \mathbf{A}$ is only positive semidefinite, meaning that $E(\mathbf{p})$ is convex but not strictly convex. In this case, there may be multiple minimizers (a subspace of solutions).

The Schoenberg--Whitney conditions~\cite{deBoor2001,Schumaker2007} provide a classical criterion for ensuring the non-singularity of the collocation matrix in B-spline interpolation. They require each collocation point to lie strictly within the support of its corresponding basis function, ensuring linear independence and full rank. While originally formulated for interpolation, these conditions are also relevant in least squares approximation, where they help ensure that the collocation matrix $\mathbf{A}$ has full column rank. In the NURBS case, this remains valid if the underlying B-spline basis satisfies the conditions and the weights are positive and well-scaled. Under these assumptions, the NURBS basis functions are independent at the data points, and $\mathbf{A}^T \mathbf{A}$ is positive definite, which guarantees a unique solution to the least squares problem.

\subsection{Geometric iterative procedure with adaptive weights}
To solve the above least squares fitting problem, an iterative procedure is use. This procedure iteratively adjusts control points and constructs a sequence of fitting patches $(\mathbf{c}^{(k)}(\mathbf{t}))_{k=0}^\infty$.

The procedure is initialized by selecting an initial set of control points $\mathbf{p}_i^{(0)}$ which are used to construct the \emph{initial fitting patch}
\begin{equation*}
    \mathbf{c}^{(0)}(\mathbf{t}) = \sum_{i \in \mathcal{I}}{\mathbf{p}_i^{(0)} B_i(\mathbf{t})}, \quad \mathbf{t} \in \Omega.
\end{equation*}
Choosing the initial control points appropriately is important to ensure rapid convergence of the iterative procedure~\cite{Deng2014}.

At each iteration, the \emph{difference vectors} between the data points and the current fitting patch $\mathbf{c}^{(k)}(\mathbf{t})$ are calculated as
\begin{equation*}
    \boldsymbol{\delta}_j^{(k)} = \mathbf{q}_j - \mathbf{c}^{(k)}(\mathbf{t}_j), \quad j \in \mathcal{J}.
\end{equation*}
These difference vectors measure the discrepancies between the data points and the fitting patch. The \emph{adjusting vectors}, which determine how the control points are updated, are computed as
\begin{equation*}
    \boldsymbol{\varDelta}_i^{(k)} = -\sum_{j \in \mathcal{J}}{B_i(\mathbf{t}_j)\boldsymbol{\delta}_j^{(k)}}, \quad i \in \mathcal{I}.
\end{equation*}
The control points are then updated as
\begin{equation*}
    \mathbf{p}_i^{(k+1)} = \mathbf{p}_i^{(k)} - \mu_i^{(k+1)}\boldsymbol{\varDelta}_i^{(k)}, \quad i \in \mathcal{I},
\end{equation*}
where $\mu_i^{(k+1)}$ are the \emph{adaptive weights} that dynamically adjust during the iterations. The sequence of fitting patches is thus generated as
\begin{equation*}
    \mathbf{c}^{(k+1)}(\mathbf{t}) = \sum_{i \in \mathcal{I}}{\mathbf{p}_i^{(k+1)} B_i(\mathbf{t})}, \quad \mathbf{t} \in \Omega.
\end{equation*}

The control point update rule can be expressed as
\begin{equation}
\label{eq:upd}
    \mathbf{p}^{(k+1)} = \mathbf{p}^{(k)} - \mathbf{W}^{(k+1)}\boldsymbol{\Delta}^{(k)},
\end{equation}
where $\mathbf{p}^{(k)} = [\mathbf{p}_i^{(k)}]_{i \in \mathcal{I}}^T$, $\mathbf{W}^{(k+1)} = \diag\left(\mu_i^{(k+1)}\right)_{i \in \mathcal{I}}$, and $\boldsymbol{\Delta}^{(k)} = [\boldsymbol{\varDelta}_i^{(k)}]_{i \in \mathcal{I}}^T$.

\citet{Rios2022}~established a connection between LSPIA and stochastic gradient descent by demonstrating that the adjusting vectors used in LSPIA correspond to components of the gradient vector of the objective function. Considering the objective function $E(\mathbf{p})$ defined by \eqref{eq:obj}, the adjusting vectors  $\boldsymbol{\varDelta}_i^{(k)}$ align with the gradient components, such that $\nabla E(\mathbf{p}^{(k)}) = \boldsymbol{\Delta}^{(k)}$.

\subsection{Adaptive weighting mechanism}
A critical enhancement in AdagradLSPIA compared to standard LSPIA is the introduction of an adaptive weighting mechanism inspired by the Adagrad optimization algorithm. Unlike standard LSPIA, which uses a fixed global weight to scale the adjusting vectors, AdagradLSPIA employs separate adaptive weights for each control point adjustment. These weights are defined as
\begin{equation}
\label{eq:aw}
    \mu^{(k+1)}_i = \frac{\mu}{\sqrt{\epsilon + v^{(k+1)}_i}}, \quad i \in \mathcal{I},
\end{equation}
where:
\setdefaultleftmargin{1em}{}{}{}{}{}
\begin{compactitem}
    \item $\mu$ is a global weight (also known as the learning rate in the context of machine learning) that controls the overall weight (step size);
    \item $\epsilon > 0$ is a small constant that ensures numerical stability, preventing division by zero;
    \item $v^{(k+1)}_i$ is the historical (accumulated) gradient information, calculated as
    \begin{equation*}      
        v^{(0)}_i = 0, \quad v^{(k+1)}_i = v^{(k)}_i + \lVert \boldsymbol{\varDelta}_i^{(k)} \rVert^2, \quad i \in \mathcal{I}.
    \end{equation*}
\end{compactitem}

The adaptive weights $\mu^{(k+1)}_i$ allow the method to take larger steps in directions where the gradients are smaller, thereby accelerating convergence in smoother regions of the fitting patch. In contrast, the weights decrease when the gradients are large, ensuring that the method remains stable and does not overshoot during the iterative process.

\subsection{Convergence}
The convergence of AdagradLSPIA can be understood as a natural extension of the well-studied convergence behavior of the Adagrad optimization algorithm to the geometric iterative framework of LSPIA. In the convex, smooth setting of least squares fitting, Adagrad is known to guarantee that the iterative sequence converges to the unique minimizer, due to its adaptive per-coordinate step sizes and the associated summability of weighted gradient norms~\cite{Duchi2011}. By integrating this mechanism into LSPIA, AdagradLSPIA inherits the descent and quasi-Fej\'{e}r monotonicity properties that underlie Adagrad's convergence guarantees, while simultaneously preserving the proven convergence of standard LSPIA in reducing the least squares error~\cite{Rios2022}. Consequently, AdagradLSPIA not only accelerates the decrease of fitting error through adaptive weighting but also retains the theoretical convergence guarantees of both Adagrad and LSPIA, making it particularly robust and efficient for convex problems (for example, tensor product B-spline surface fitting tasks).

\subsection{Hyperparameter tuning}
The performance of the AdagradLSPIA method depends on the appropriate selection of hyperparameters, particularly $\mu$ and $\epsilon$. The global weight $\mu$ determines the scale of adjustments to the control points, and it is important to balance this parameter to avoid overly aggressive updates, which can lead to oscillations, or overly conservative updates, which may slow down convergence. In practice, $\mu$ can be selected by cross-validation or empirical testing.

The stability parameter $\epsilon$ is usually set to a small value (e.g. $10^{-8}$) to avoid division by zero in the calculation of $\mu^{(k+1)}_i$. Although $\epsilon$ does not have a significant impact on convergence behavior, it can be adjusted based on the scale of the gradients to improve numerical stability.

\subsection{Computational complexity}

The computational complexity of the AdagradLSPIA method is determined by its iterative nature and adaptive weighting mechanism. In each iteration, the key computational steps involve evaluating the points on the fitting patch, computing the difference vectors, determining the adjustment vectors, and updating the control points using the adjustment vectors and adaptive weights. The complexity of these operations depends on the number of data points $\lvert \mathcal{J} \rvert$, the number of control points $\lvert \mathcal{I} \rvert$, and the computational cost associated with the basis function evaluation.

In the standard LSPIA method, the computational cost of each iteration is dominated by the calculation of adjusting vectors, which requires summing over all data points for each control point. Therefore, the computational complexity is $\mathcal{O}(\lvert \mathcal{I} \rvert \cdot \lvert \mathcal{J} \rvert)$ per iteration. AdagradLSPIA introduces an additional step of computing and updating adaptive weights, which involves maintaining and scaling the historical gradient information. Although this adds a minor computational overhead, the total complexity per iteration remains $\mathcal{O}(\lvert \mathcal{I} \rvert \cdot \lvert \mathcal{J} \rvert)$, as the extra operations scale linearly with $\lvert \mathcal{I} \rvert$ and $\lvert \mathcal{J} \rvert$.

When comparing AdagradLSPIA with other enhanced LSPIA variants, such as momentum-based methods (e.g., PmLSPIA and NmLSPIA)~\cite{Liu2024} or memory-augmented LSPIA (MLSPIA)~\cite{Huang2020}, the computational complexity remains competitive. Momentum-based methods involve additional computations to store and update momentum terms, while MLSPIA requires maintaining a history of previous adjustments. These methods typically incur similar or slightly higher costs than AdagradLSPIA due to their additional storage and update requirements. However, dynamic scaling in AdagradLSPIA offers a unique advantage in accelerating convergence, reducing the total number of iterations required and thereby compensating for the cost per iteration.

High-order iterative methods, such as hyperpower LSPIA (HPLSPIA)~\cite{Sajavicius2023}, achieve faster convergence by leveraging advanced matrix computations, but these come with significantly higher per-iteration costs, often exceeding $\mathcal{O}(\lvert \mathcal{I} \rvert^2 \cdot \lvert \mathcal{J} \rvert)$. In contrast, AdagradLSPIA strikes a balance between computational efficiency and convergence rate by maintaining linear complexity per iteration while achieving superior stability and accuracy through its adaptive weighting mechanism.

\section{Experimental study}
\label{sec:3}

This section details the implementation of the AdagradLSPIA method and presents the experimental results of its comparison with the standard LSPIA method in the context of tensor product B-spline surface fitting. The experiments were conducted on both noise-free and noisy data to examine the performance and robustness of the proposed method.  

\subsection{Implementation details}
\begin{figure*}
\centering
\includegraphics[scale=1.0]{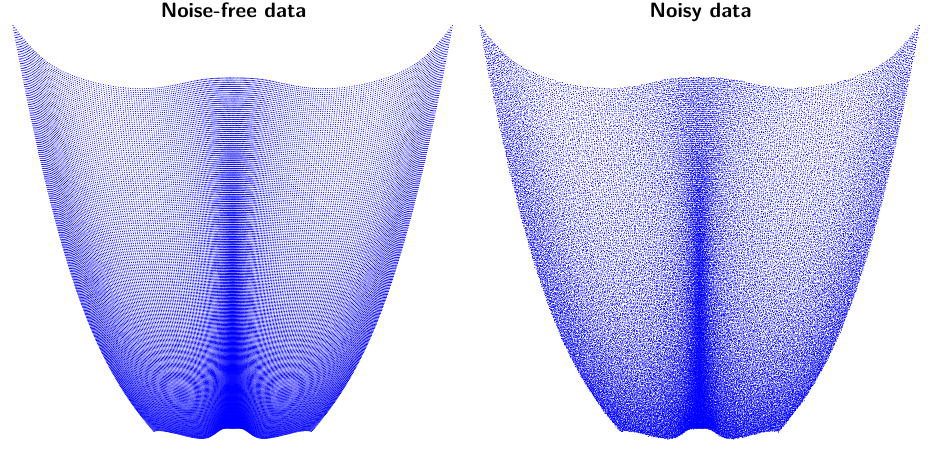}
\caption{Data points representing an automobile hood model, including noise-free data (left) and noisy data (right) with added zero-mean Gaussian noise (standard deviation 0.02).\label{fig:01}}
\end{figure*}

The experiments utilize a topologically rectangular set of $251 \times 251$ data points representing an automobile hood that features a combination of smooth regions, sharp transitions, and subtle design details~(see~\autoref{fig:01}). For the noisy data case, zero-mean Gaussian noise with a standard deviation of $0.02$ was added to the original data. Parameter values for the data points were assigned using a uniform parameterization method in both parametric directions.

A bicubic tensor product B-spline surface with $25 \times 25$ control points was used for the fitting. The knot vectors to define the B-spline basis functions were constructed using the averaging technique (see~\cite{Piegl1997,Kineri2012}). For initialization, all control points were set to zero.

The fitting error was computed at each iteration as the sum of the squared Euclidean distances between the data points and their corresponding points on the fitting surface. The iterative process terminated either when the absolute difference between consecutive fitting errors was less than the prescribed tolerance ($10^{-7}$ for noise-free data and $10^{-7}$ for noisy data) or when the maximum of $10^3$ iterations was reached.

\subsection{Results and discussion}

The experiments analyzed the performance of LSPIA and AdagradLSPIA on noise-free and noisy data by examining their sensitivity to the global weight, convergence behavior, and final surface quality.

\begin{figure*}
\centering
\vskip -72pt
\includegraphics[scale=1.0]{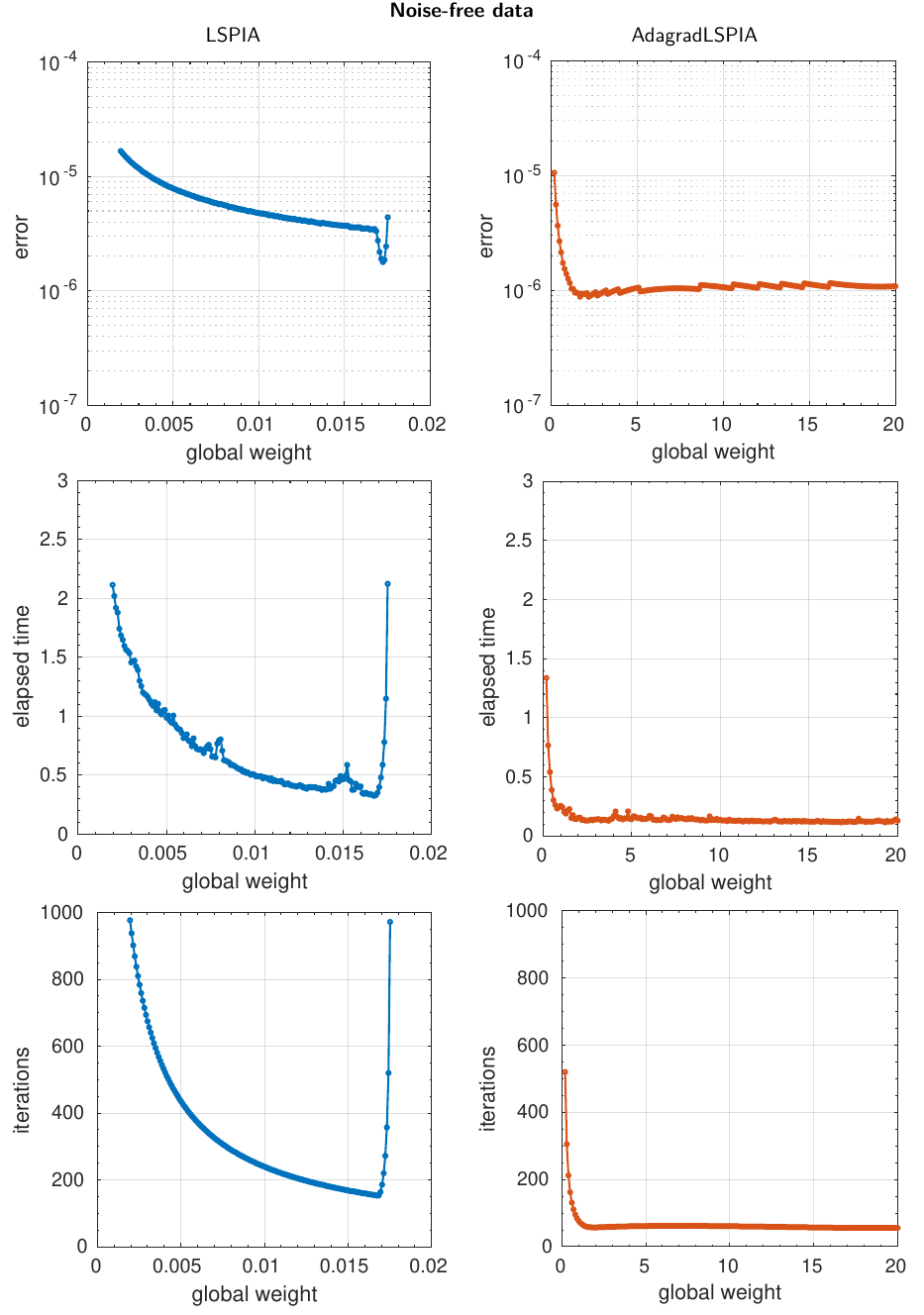}
\vskip -12pt
\caption{Comparison of the LSPIA and AdagradLSPIA methods for noise-free data across a range of global weight values. The plots depict the variations in fitting error, elapsed time (in s), and the number of iterations required to converge, highlighting the broader convergence region and superior performance of AdagradLSPIA.\label{fig:02}}
\end{figure*}

\begin{figure*}
\centering
\vskip -72pt
\includegraphics[scale=1.0]{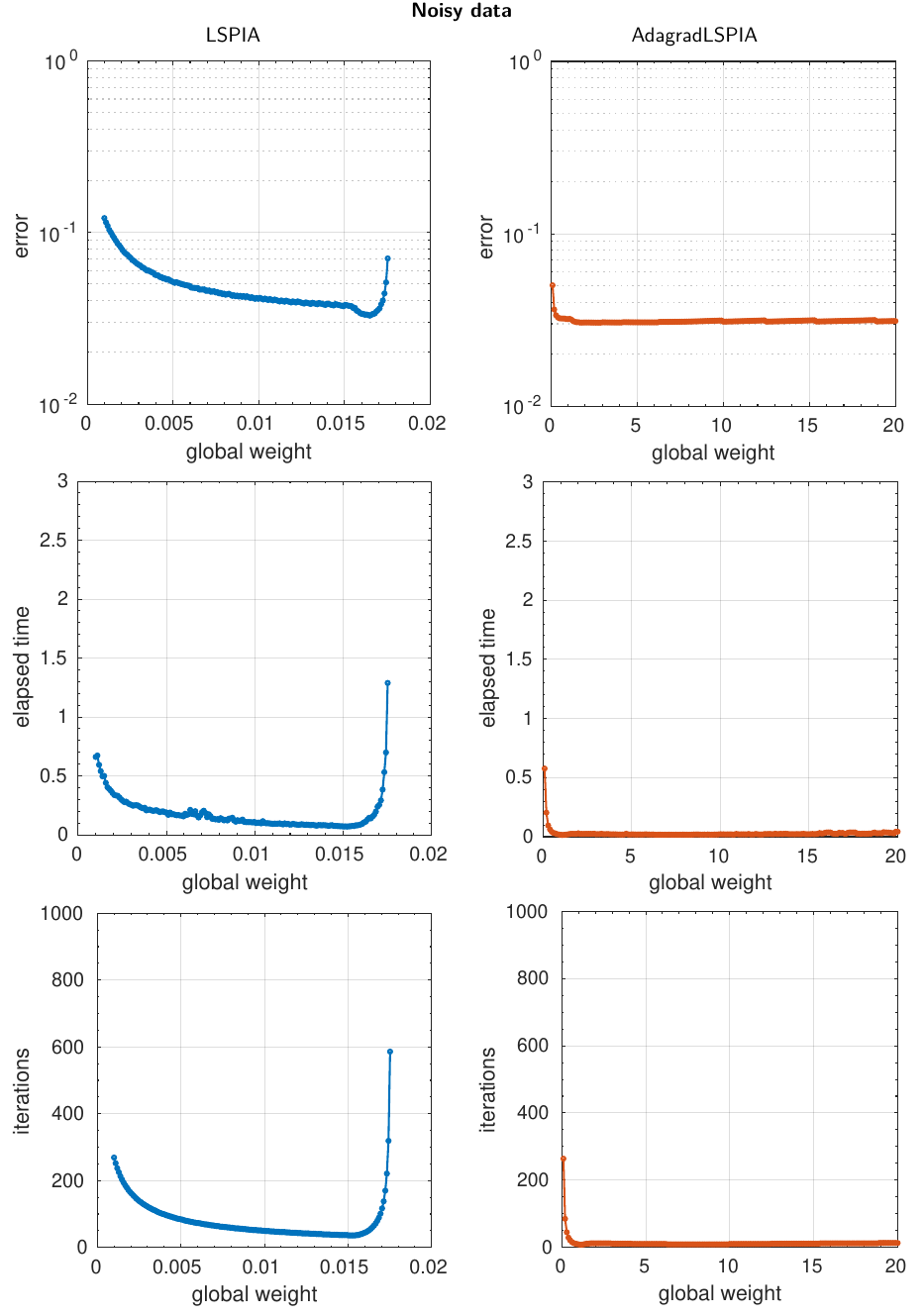}
\vskip -12pt
\caption{Comparison of the LSPIA and AdagradLSPIA methods for noisy data across a range of global weight values. The plots depict the variations in fitting error, elapsed time (in s), and the number of iterations required to converge, highlighting the broader convergence region and superior performance of AdagradLSPIA.\label{fig:03}}
\end{figure*}

\paragraph{Influence of global weight}
The initial set of experiments investigated how the global weight affected the performance of the LSPIA and AdagradLSPIA methods. LSPIA was tested with 200 global weight values sampled uniformly from the interval $(0, 0.02]$, while AdagradLSPIA was examined over $(0, 20]$. The results~(\autoref{fig:02} and \autoref{fig:03}) show that LSPIA converged only within a limited range of global weights; outside this range, it reached the maximum iteration limit without achieving the desired accuracy. In contrast, AdagradLSPIA converged consistently throughout the tested range, demonstrating its robustness to variations in global weight.

\paragraph{Optimized global weights}
From the experiments, global weight values yielding minimal fitting error, minimal elapsed time, and minimal iterations were identified for both methods. \autoref{tab:1}~summarizes these results, showing that AdagradLSPIA achieved better accuracy, reduced computational time, and required fewer iterations compared to the LSPIA method. For example, the minimal fitting error for AdagradLSPIA applied to noise-free data was $8.75530 \times 10^{-7}$, which is significantly lower than $1.78478 \times 10^{-6}$ achieved by LSPIA.

\begin{figure*}
\centering
\includegraphics[scale=1.0]{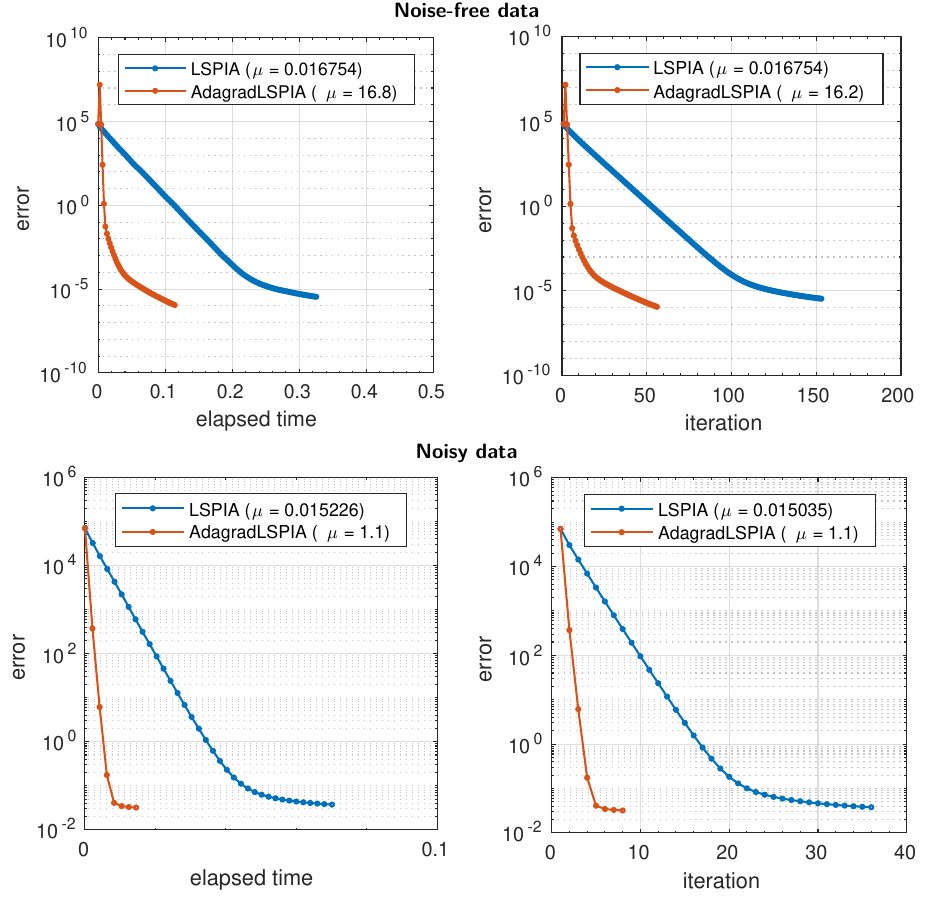}
\caption{Convergence plots for the LSPIA and AdagradLSPIA methods using optimized global weight values. The plots show the reduction in fitting error over elapsed time and iterations for both noise-free and noisy data, demonstrating the faster convergence of AdagradLSPIA. \coltext\label{fig:04}}
\end{figure*}

\paragraph{Convergence behavior}
Convergence plots~(\autoref{fig:04}) for the optimized global weights reveal sublinear convergence for both methods, with AdagradLSPIA consistently outperforming LSPIA in reducing the fitting error per iteration and per unit of elapsed time. AdagradLSPIA demonstrated faster error reduction, particularly in the initial stages, reflecting the benefits of its adaptive weighting mechanism.

\begin{figure*}
\centering
\includegraphics[scale=1.0]{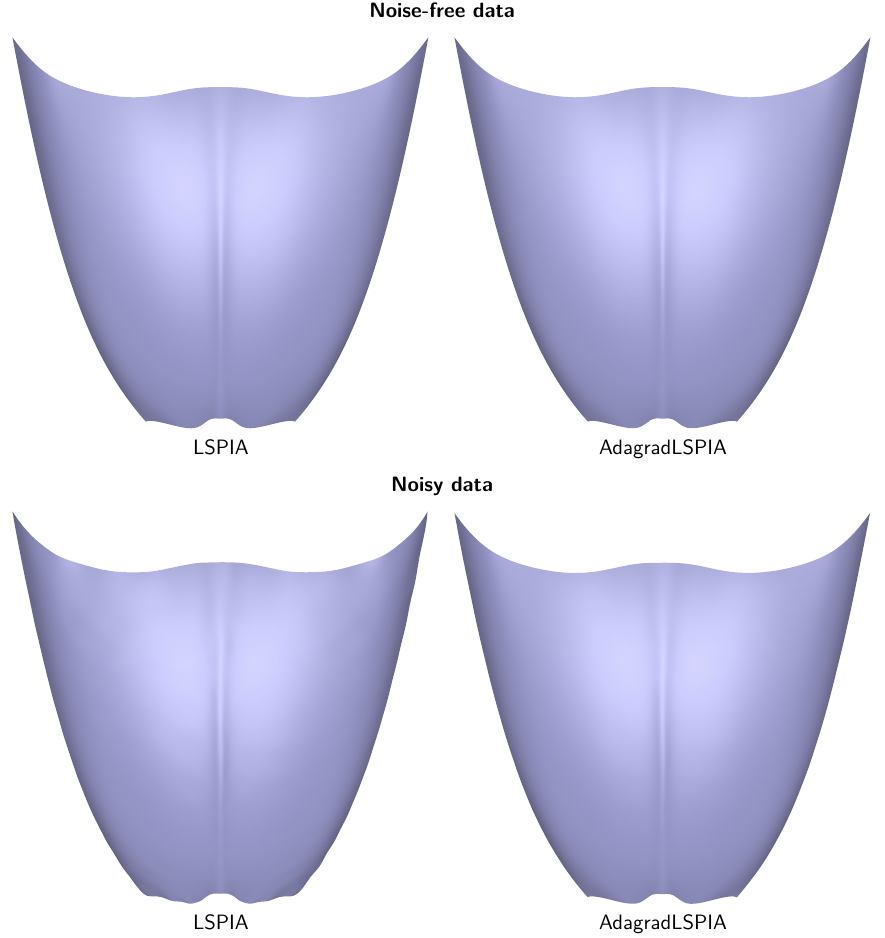}
\caption{The final fitting surfaces generated by the LSPIA and AdagradLSPIA methods for noise-free and noisy data. The surfaces were obtained using global weights that minimize the fitting error~(\autoref{tab:1}).\label{fig:05}}
\end{figure*}

\begin{figure*}
\centering
\includegraphics[scale=1.0]{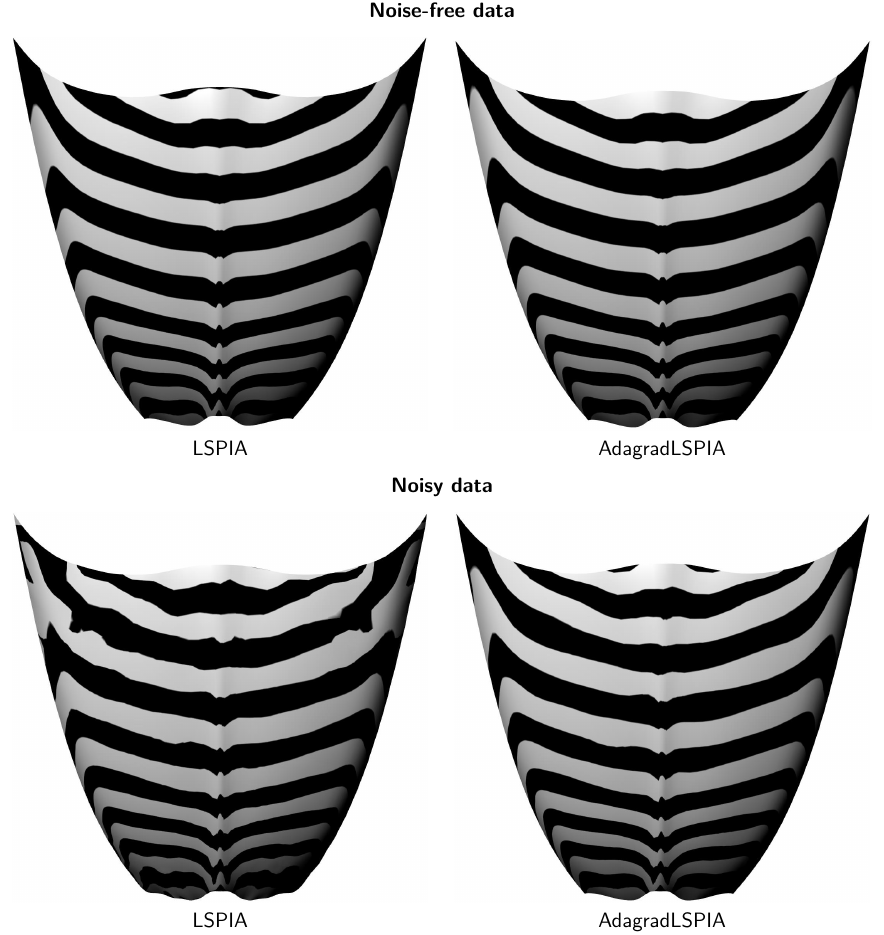}
\caption{Reflection lines (zebra mappings) on the fitting surfaces generated by the LSPIA and AdagradLSPIA methods for noise-free and noisy data. These lines demonstrate the ability of AdagradLSPIA to capture geometric details and maintain surface smoothness, even on noisy data.\label{fig:06}}
\end{figure*}

\paragraph{Surface quality}
The final fitting surfaces~(\autoref{fig:05}), obtained using the global weights corresponding to the minimal fitting errors, and their reflection lines~(\autoref{fig:06}) provide a qualitative comparison of the results. The reflection lines show the ability of AdagradLSPIA to capture geometric details and maintain surface smoothnes, even on noisy data.

\begin{table}
\centering
\caption{\label{tab:1}  Performance comparison of the LSPIA and AdagradLSPIA methods for the automobile hood model. The table includes the initial fitting error, final fitting errors, elapsed times (in s), and the number of iterations required for global weights corresponding to minimal fitting error, elapsed time, and iterations.}
\vskip 5pt
\footnotesize
\begin{tabular}{p{2.5cm} p{2.5cm} p{2.5cm} p{2.5cm} p{2.5cm}}
\hline
Method & Global weight & Fitting error & Elapsed time & Iterations  \\
\hline
\multicolumn{5}{l}{\emph{Initial fitting surface}} \\
 & & 7.03910E+04 & & \\[3pt]
\multicolumn{5}{l}{\textbf{Noise-free data}} \\
\multicolumn{5}{l}{\emph{Minimal fitting error}} \\
LSPIA & 0.017231 & 1.78478E-06 & 0.588365 & 272 \\
AdagradLSPIA & 2.200000 & 8.75530E-07 & 0.135451 & 58 \\[3pt]
\multicolumn{5}{l}{\emph{Minimal elapsed time}} \\
LSPIA & 0.016754 & 3.47095E-06 & 0.325270 & 153 \\
AdagradLSPIA & 16.800000 & 1.12510E-06 & 0.114263 & 56 \\[3pt]
\multicolumn{5}{l}{\emph{Minimal number of iterations}} \\
LSPIA & 0.016754 & 3.47095E-06 & 0.325270 & 153 \\
AdagradLSPIA & 16.200000 & 1.15355E-06 & 0.122095 & 56 \\[3pt]
\multicolumn{5}{l}{\textbf{Noiy data}} \\
\multicolumn{5}{l}{\emph{Minimal fitting error}} \\
LSPIA & 0.016467 & 3.27684E-02 & 0.144296 & 61 \\
AdagradLSPIA & 2.900000 & 3.04729E-02 & 0.025656 & 12 \\[3pt]
\multicolumn{5}{l}{\emph{Minimal elapsed time}} \\
LSPIA & 0.015226 & 3.74664E-02 & 0.070095 & 36 \\
AdagradLSPIA & 1.100000 & 3.20504E-02 & 0.014523 & 8 \\[3pt]
\multicolumn{5}{l}{\emph{Minimal number of iterations}} \\
LSPIA & 0.015035 & 3.78637E-02 & 0.070474 & 36 \\
AdagradLSPIA & 1.100000 & 3.20504E-02 & 0.014523 & 8 \\
\hline
\end{tabular}
\end{table}

\subsection{Summary of results}
The experiments demonstrate the advantages of AdagradLSPIA over the standard LSPIA method. AdagradLSPIA consistently achieved higher accuracy, faster convergence, and improved surface quality. Its robustness to global weight variations and adaptability to the geometric characteristics of the data make it a superior choice for surface fitting tasks.

\section{Concluding remarks}
\label{sec:4}

This paper introduced the Adaptive Gradient Least Squares Progressive Iterative Approximation (AdagradLSPIA) method, a significant enhancement of the standard LSPIA technique. By integrating the Adagrad optimization algorithm into the geometric iterative framework, AdagradLSPIA dynamically adjusts weights based on historical gradient information, resulting in more stable and efficient convergence.

The experimental results confirmed the practical advantages of the AdagradLSPIA method, including faster convergence, reduced computational cost, and improved accuracy compared to the standard approach. In particular, AdagradLSPIA exhibited superior performance in a range of global weight values and consistently produced smooth and accurate fitting surfaces, as validated by reflection line analysis.

The computational complexity of AdagradLSPIA remains consistent with the standard LSPIA method, with the adaptive weighting mechanism introducing minimal additional computational overhead. Compared to other enhanced LSPIA variants, it offers an efficient and practical solution for surface fitting tasks, especially when faster convergence and robustness are required.

In addition to its promising results, AdagradLSPIA opens avenues for further research. The integration of other advanced adaptive optimizers, such as Adam~\cite{Kingma2014}, Adadelta~\cite{Zeiler2012}, Nadam~\cite{Dozat2016}, hypergradient descent~\cite{Baydin2018}, or AMSGrad~\cite{Reddi2018}, could potentially improve the efficiency and robustness of the method. Future work may also extend the method to handle more complex geometries.

Overall, AdagradLSPIA represents a robust and efficient tool for surface fitting, with significant potential for applications in geometric modeling, CAD, and reverse engineering.

\bibliography{references_iga,references_ml}

\begin{thebibliography}{25}
\expandafter\ifx\csname natexlab\endcsname\relax\def\natexlab#1{#1}\fi
\providecommand{\url}[1]{\texttt{#1}}
\providecommand{\href}[2]{#2}
\providecommand{\path}[1]{#1}
\providecommand{\DOIprefix}{doi:}
\providecommand{\ArXivprefix}{arXiv:}
\providecommand{\URLprefix}{URL: }
\providecommand{\Pubmedprefix}{pmid:}
\providecommand{\doi}[1]{\href{http://dx.doi.org/#1}{\path{#1}}}
\providecommand{\Pubmed}[1]{\href{pmid:#1}{\path{#1}}}
\providecommand{\bibinfo}[2]{#2}
\ifx\xfnm\relax \def\xfnm[#1]{\unskip,\space#1}\fi
\bibitem[{Deng and Lin(2014)}]{Deng2014}
\bibinfo{author}{C.~Deng}, \bibinfo{author}{H.~Lin},
\newblock \bibinfo{title}{Progressive and iterative approximation for least
  squares {B}-spline curve and surface fitting},
\newblock \bibinfo{journal}{Comput.-Aided Des.} \bibinfo{volume}{47}
  (\bibinfo{year}{2014}) \bibinfo{pages}{32--44}.
\bibitem[{Lin et~al.(2004)Lin, Wang, and Dong}]{Lin2004}
\bibinfo{author}{H.~Lin}, \bibinfo{author}{G.~Wang}, \bibinfo{author}{C.~Dong},
\newblock \bibinfo{title}{Constructing iterative non-uniform {B}-spline curve
  and surface to fit data points},
\newblock \bibinfo{journal}{Sci. China Ser. F} \bibinfo{volume}{47}
  (\bibinfo{year}{2004}) \bibinfo{pages}{315--331}.
\bibitem[{Lin et~al.(2018)Lin, Maekawa, and Deng}]{Lin2018}
\bibinfo{author}{H.~Lin}, \bibinfo{author}{T.~Maekawa},
  \bibinfo{author}{C.~Deng},
\newblock \bibinfo{title}{Survey on geometric iterative methods and their
  applications},
\newblock \bibinfo{journal}{Comput.-Aided Des.} \bibinfo{volume}{95}
  (\bibinfo{year}{2018}) \bibinfo{pages}{40--51}.
\bibitem[{Huang and Wang(2020)}]{Huang2020}
\bibinfo{author}{Z.-D. Huang}, \bibinfo{author}{H.-D. Wang},
\newblock \bibinfo{title}{On a progressive and iterative approximation method
  with memory for least square fitting},
\newblock \bibinfo{journal}{Comput. Aided Geom. Design} \bibinfo{volume}{82}
  (\bibinfo{year}{2020}) \bibinfo{pages}{101931}.
\bibitem[{Wu and Liu(2024)}]{Wu2024}
\bibinfo{author}{N.-C. Wu}, \bibinfo{author}{C.~Liu},
\newblock \bibinfo{title}{Asynchronous progressive iterative approximation
  method for least squares fitting},
\newblock \bibinfo{journal}{Comput. Aided Geom. Design} \bibinfo{volume}{111}
  (\bibinfo{year}{2024}) \bibinfo{pages}{102295}.
\bibitem[{Liu et~al.(2024)Liu, Wu, Li, and Hu}]{Liu2024}
\bibinfo{author}{C.~Liu}, \bibinfo{author}{N.-C. Wu}, \bibinfo{author}{J.~Li},
  \bibinfo{author}{L.~Hu},
\newblock \bibinfo{title}{Two novel iterative approaches for improved {LSPIA}
  convergence},
\newblock \bibinfo{journal}{Comput. Aided Geom. Design} \bibinfo{volume}{111}
  (\bibinfo{year}{2024}) \bibinfo{pages}{102312}.
\bibitem[{Song and Bo(2024)}]{Song2024}
\bibinfo{author}{Q.~Song}, \bibinfo{author}{P.~Bo},
\newblock \bibinfo{title}{Newton geometric iterative method for {B}-spline
  curve and surface approximation},
\newblock \bibinfo{journal}{Comput.-Aided Des.} \bibinfo{volume}{172}
  (\bibinfo{year}{2024}) \bibinfo{pages}{103716}.
\bibitem[{Ebrahimi and Loghmani(2019)}]{Ebrahimi2019}
\bibinfo{author}{A.~Ebrahimi}, \bibinfo{author}{G.~B. Loghmani},
\newblock \bibinfo{title}{A composite iterative procedure with fast convergence
  rate for the progressive-iteration approximation of curves},
\newblock \bibinfo{journal}{J. Comput. Appl. Math.} \bibinfo{volume}{359}
  (\bibinfo{year}{2019}) \bibinfo{pages}{1--15}.
\bibitem[{Sajavi\v{c}ius(2023)}]{Sajavicius2023}
\bibinfo{author}{S.~Sajavi\v{c}ius},
\newblock \bibinfo{title}{Hyperpower least squares progressive iterative
  approximation},
\newblock \bibinfo{journal}{J. Comput. Appl. Math.} \bibinfo{volume}{422}
  (\bibinfo{year}{2023}) \bibinfo{pages}{114888}.
\bibitem[{Yao and Hu(2025)}]{Yao2025}
\bibinfo{author}{Z.~Yao}, \bibinfo{author}{Q.~Hu},
\newblock \bibinfo{title}{Accelerated local progressive-iterative approximation
  methods for curve and surface fitting},
\newblock \bibinfo{journal}{The Visual Computer}  (\bibinfo{year}{2025}).
\bibitem[{Yao and Hu(2024)}]{Yao2024}
\bibinfo{author}{Z.~Yao}, \bibinfo{author}{Q.~Hu},
\newblock \bibinfo{title}{Distributed least-squares progressive iterative
  approximation for blending curves and surfaces},
\newblock \bibinfo{journal}{Comput.-Aided Des.} \bibinfo{volume}{175}
  (\bibinfo{year}{2024}) \bibinfo{pages}{103749}.
\bibitem[{Lan et~al.(2024)Lan, Ji, Wang, and Zhu}]{Lan2024}
\bibinfo{author}{L.~Lan}, \bibinfo{author}{Y.~Ji}, \bibinfo{author}{M.-Y.
  Wang}, \bibinfo{author}{C.-G. Zhu},
\newblock \bibinfo{title}{{Full-LSPIA}: {A} least-squares progressive-iterative
  approximation method with optimization of weights and knots for nurbs curves
  and surfaces},
\newblock \bibinfo{journal}{Comput.-Aided Des.} \bibinfo{volume}{169}
  (\bibinfo{year}{2024}) \bibinfo{pages}{103673}.
\bibitem[{Duchi et~al.(2011)Duchi, Hazan, and Singer}]{Duchi2011}
\bibinfo{author}{J.~Duchi}, \bibinfo{author}{E.~Hazan},
  \bibinfo{author}{Y.~Singer},
\newblock \bibinfo{title}{Adaptive subgradient methods for online learning and
  stochastic optimization},
\newblock \bibinfo{journal}{J. Mach. Learn. Res.} \bibinfo{volume}{12}
  (\bibinfo{year}{2011}) \bibinfo{pages}{2121--2159}.
\bibitem[{Bottou et~al.(2018)Bottou, Curtis, and Nocedal}]{Bottou2018}
\bibinfo{author}{L.~Bottou}, \bibinfo{author}{F.~E. Curtis},
  \bibinfo{author}{J.~Nocedal},
\newblock \bibinfo{title}{Optimization methods for large-scale machine
  learning},
\newblock \bibinfo{journal}{SIAM Rev.} \bibinfo{volume}{60}
  (\bibinfo{year}{2018}) \bibinfo{pages}{223--311}.
\bibitem[{Kochenderfer and Wheeler(2019)}]{Kochenderfer2019}
\bibinfo{author}{M.~J. Kochenderfer}, \bibinfo{author}{T.~A. Wheeler},
  \bibinfo{title}{Algorithms for Optimization}, \bibinfo{publisher}{The MIT
  Press}, \bibinfo{year}{2019}.
\bibitem[{de~{B}oor(2001)}]{deBoor2001}
\bibinfo{author}{C.~de~{B}oor}, \bibinfo{title}{A Practical Guide to Splines},
  \bibinfo{edition}{revised} ed., \bibinfo{publisher}{Springer-Verlag New
  York}, \bibinfo{year}{2001}.
\bibitem[{Schumaker(2007)}]{Schumaker2007}
\bibinfo{author}{L.~Schumaker}, \bibinfo{title}{Spline Functions: Basic
  Theory}, \bibinfo{edition}{3rd} ed., \bibinfo{publisher}{Cambridge University
  Press}, \bibinfo{year}{2007}.
\bibitem[{Rios and J\"uttler(2022)}]{Rios2022}
\bibinfo{author}{D.~Rios}, \bibinfo{author}{B.~J\"uttler},
\newblock \bibinfo{title}{{LSPIA}, (stochastic) gradient descent, and parameter
  correction},
\newblock \bibinfo{journal}{J. Comput. Appl. Math.} \bibinfo{volume}{406}
  (\bibinfo{year}{2022}) \bibinfo{pages}{113921}.
\bibitem[{Piegl and Tiller(1997)}]{Piegl1997}
\bibinfo{author}{L.~Piegl}, \bibinfo{author}{W.~Tiller}, \bibinfo{title}{The
  {NURBS} Book}, Monographs in Visual Communications, \bibinfo{edition}{2nd}
  ed., \bibinfo{publisher}{Springer Berlin Heidelberg}, \bibinfo{year}{1997}.
\bibitem[{Kineri et~al.(2012)Kineri, Wang, Lin, and Maekawa}]{Kineri2012}
\bibinfo{author}{Y.~Kineri}, \bibinfo{author}{M.~Wang},
  \bibinfo{author}{H.~Lin}, \bibinfo{author}{T.~Maekawa},
\newblock \bibinfo{title}{${B}$-spline surface fitting by iterative geometric
  interpolation/approximation algorithms},
\newblock \bibinfo{journal}{Comput.-Aided Des.} \bibinfo{volume}{44}
  (\bibinfo{year}{2012}) \bibinfo{pages}{697--708}.
\bibitem[{Kingma and Ba(2014)}]{Kingma2014}
\bibinfo{author}{D.~P. Kingma}, \bibinfo{author}{J.~Ba}, \bibinfo{title}{Adam:
  A method for stochastic optimization}, \bibinfo{year}{2014}. \URLprefix
  \url{https://arxiv.org/abs/1412.6980}.
  \DOIprefix\doi{10.48550/ARXIV.1412.6980}.
\bibitem[{Zeiler(2012)}]{Zeiler2012}
\bibinfo{author}{M.~D. Zeiler}, \bibinfo{title}{{ADADELTA}: An adaptive
  learning rate method}, \bibinfo{year}{2012}. \URLprefix
  \url{https://arxiv.org/abs/1212.5701}.
  \DOIprefix\doi{10.48550/ARXIV.1212.5701}.
\bibitem[{Dozat(2016)}]{Dozat2016}
\bibinfo{author}{T.~Dozat},
\newblock \bibinfo{title}{Incorporating {N}esterov momentum into {A}dam},
\newblock in: \bibinfo{booktitle}{International Conference on Learning
  Representations}, \bibinfo{year}{2016}. \URLprefix
  \url{https://openreview.net/forum?id=OM0jvwB8jIp57ZJjtNEZ}.
\bibitem[{Baydin et~al.(2018)Baydin, Cornish, Rubio, Schmidt, and
  Wood}]{Baydin2018}
\bibinfo{author}{A.~G. Baydin}, \bibinfo{author}{R.~Cornish},
  \bibinfo{author}{D.~M. Rubio}, \bibinfo{author}{M.~Schmidt},
  \bibinfo{author}{F.~Wood},
\newblock \bibinfo{title}{Online learning rate adaptation with hypergradient
  descent},
\newblock in: \bibinfo{booktitle}{International Conference on Learning
  Representations}, \bibinfo{year}{2018}. \URLprefix
  \url{https://openreview.net/forum?id=BkrsAzWAb}.
\bibitem[{Reddi et~al.(2018)Reddi, Kale, and Kumar}]{Reddi2018}
\bibinfo{author}{S.~J. Reddi}, \bibinfo{author}{S.~Kale},
  \bibinfo{author}{S.~Kumar},
\newblock \bibinfo{title}{On the convergence of {A}dam and beyond},
\newblock in: \bibinfo{booktitle}{International Conference on Learning
  Representations}, \bibinfo{year}{2018}. \URLprefix
  \url{https://openreview.net/forum?id=ryQu7f-RZ}.

\end{thebibliography}
\addcontentsline{toc}{section}{References}

\end{document}